\documentclass[12pt]{article}%
\usepackage{amsfonts}
\usepackage{sw20bams}
\usepackage{amsmath}
\usepackage{amssymb}
\usepackage{graphicx}%
\providecommand{\U}[1]{\protect\rule{.1in}{.1in}}
\begin{document}

\title{Abel's Functional Equation and Interrelations}
\author{Steven Finch}
\date{March 18, 2025}
\maketitle

\begin{abstract}
Convex solutions $A,B,I,J$ of four Abel equations are numerically studied.
\ We do not know exact formulas for any of these functions, but conjecture
that $A,B$ and $I,J$ are closely related. \ [Corrigendum at end.]

\end{abstract}

\footnotetext{Copyright \copyright \ 2025 by Steven R. Finch. All rights
reserved.}

Let $f:(0,r)\rightarrow(0,\infty)$ be concave, where $0<r\leq\infty$; assume
that $f$ obeys%
\[%
\begin{array}
[c]{ccccc}%
0<f(x)<x &  & \text{always and} &  & \lim\limits_{x\rightarrow0^{+}}%
\dfrac{f(x)}{x}=1.
\end{array}
\]
The Abel equation%
\[
F\left(  f(x)\right)  =F(x)+1
\]
hence possesses a convex solution $F:(0,r)\rightarrow(-\infty,\infty)$, unique
up to an additive constant. \ If $f$ is strictly increasing, then $F$ is
strictly decreasing \cite{K1-abel, K2-abel, K3-abel}. \ If $f$ is real
analytic, then $F$ is real analytic as well \cite{BD-abel}. \ 

For example,%
\[%
\begin{array}
[c]{ccccc}%
F\left(  \dfrac{x}{1+x}\right)  =F(x)+1 &  & \text{implies} &  &
F(x)=\dfrac{1}{x}%
\end{array}
\]
for the domain $(0,\infty)$. \ We omit the additive constant in our expression
for $F$, in the interest of simplicity. \ More generally, given $s>0$,%
\[%
\begin{array}
[c]{ccccc}%
F\left(  \dfrac{x}{(1+x^{s})^{1/s}}\right)  =F(x)+1 &  & \text{implies} &  &
F(x)=\dfrac{1}{x^{s}}.
\end{array}
\]
Closed-form solutions as such seem to be exceedingly rare. \ The special case
$s=1$ plays a featured role in the ensuing. \ Also define $\varphi=\left(
1+\sqrt{5}\right)  /2=1.618...$ (the Golden mean) and $\psi=0.629...$ (which
satisfies $I^{\prime\prime}(\psi)=0$, to be clarified shortly).

Consider now the convex solutions of four Abel equations:%
\[%
\begin{array}
[c]{ccccc}%
A\left(  x(1-x)\right)  =A(x)+1 &  & \text{where }r=1 &  & \text{(}A\text{ has
a minimum at }x=1/2\text{);}%
\end{array}
\]%
\[%
\begin{array}
[c]{ccccc}%
B\left(  \dfrac{x}{1+x+x^{2}}\right)  =B(x)+1 &  & \text{where }r=\infty &  &
\text{(}B\text{ has a minimum at }x=1\text{);}%
\end{array}
\]%
\[%
\begin{array}
[c]{ccccc}%
I\left(  \dfrac{x}{1+x-x^{2}}\right)  =I(x)+1 &  & \text{where }r=\psi &  &
\text{(}I\text{ has an inflection point at }x=\psi\text{);}%
\end{array}
\]%
\[%
\begin{array}
[c]{ccccc}%
J\left(  \dfrac{x(1+x)}{1+2x}\right)  =J(x)+1 &  & \text{where }r=\infty &  &
\text{(}J\text{ has no critical points).}%
\end{array}
\]
While we cannot determine exact expressions for $A$, $B$, $I$ or $J$, we might
try at least to find interrelations between these. \ As far as is known, the
following two conjectured identities are new:%
\[%
\begin{array}
[c]{ccc}%
A\left(  \dfrac{x}{1+x}\right)  =B(x)+1, &  & I\left(  \dfrac{x}{1+x}\right)
=J(x)+1
\end{array}
\]
over the domain $(0,\infty)$. \ The evidence underlying these surprising
formulas is entirely numerical and shall occupy us in Section 1. \ It confirms
the existence and convexity of $A$ for $0<x<1$, $B$ and $J$ for $0<x<\infty$,
and $I$ for $0<x<\psi$. \ It supports the existence and concavity of $I$ for
$\psi<x<1$ in particular; the data is less clear when $1<x<\varphi$.
\ Commonalities\ of a number-theoretic nature between $A$ \& $B$ are given in
Section 3 and likewise for $I$ \&\ $J$ in Section 4. \ These are fairly subtle
and do not immediately suggest any explanations. \ Variations of $B$ \&\ $I$
are presented in Section 5 that are noticeably distinct from the originals,
but nevertheless intractable.

\section{Iterational Asymptotics}

The logistic recurrence%
\[%
\begin{array}
[c]{ccccc}%
x_{n}=x_{n-1}\left(  1-x_{n-1}\right)  &  & \text{for }n\geq1, &  & 0<x_{0}<1
\end{array}
\]
has been analyzed via both a brute-force matching-coefficient method
\cite{Sc-abel, F1-abel} and the relatively elegant Mavecha-Laohakosol
algorithm \cite{ML-abel, F2-abel} based on work by de Bruijn \cite{dB-abel}
and Bencherif \&\ Robin \cite{BR-abel}. \ A power-logarithmic series%
\[
x_{n}\sim%
{\displaystyle\sum\limits_{m=0}^{k-1}}
P_{m}\left(  -\ln(n)-C\right)  \frac{1}{n^{m+1}}%
\]
is the asymptotic outcome of both computations, valid as $n\rightarrow\infty$,
where $P_{m}=P_{m}(X)$ is a polynomial with rational coefficients:
\[%
\begin{array}
[c]{ccccccc}%
P_{0}=1, &  & P_{1}=X, &  & P_{2}=\dfrac{1}{2}+X+X^{2}, &  & P_{3}=\dfrac
{5}{6}+\dfrac{5}{2}X+\dfrac{5}{2}X^{2}+X^{3},
\end{array}
\]%
\[%
\begin{array}
[c]{ccc}%
P_{4}=\dfrac{61}{36}+\dfrac{35}{6}X+\dfrac{15}{2}X^{2}+\dfrac{13}{3}%
X^{3}+X^{4}, &  & P_{5}=\dfrac{2609}{720}+\dfrac{515}{36}X+\dfrac{265}%
{12}X^{2}+\dfrac{101}{6}X^{3}+\dfrac{77}{12}X^{4}+X^{5},
\end{array}
\]%
\[
P_{6}=\dfrac{29069}{3600}+\dfrac{12977}{360}X+65X^{2}+61X^{3}+\dfrac{95}%
{3}X^{4}+\dfrac{87}{10}X^{5}+X^{6}%
\]
and $C=C(x_{0})$ is a constant. \ (Parameters $\tau$ and $\lambda$ in
\cite{F2-abel} are each equal to $1$ here; what was called $b_{1}$ is $1$ for
$A$ and $B$.) \ The parameter $k$ was fixed to be $4$ in \cite{F1-abel} and
$7$ in \cite{F2-abel}; here we fix $k$ to be $20$. \ 

Our simple procedure for estimating $C$, given $x_{0}$, involves computing
$x_{N}$ exactly via recursion, for some suitably large index $N$. We then set
the value $x_{N}$ equal to our series and numerically solve for $C$. \ The
assignment $x_{0}\mapsto C(x_{0})$, in the context of iterations, is written
as $x\mapsto F(x)$ when speaking of functional equations. \ Take the logistic
case, for example:%
\begin{align*}
A\left(  \frac{1}{2}\right)   &
=1.76799378613615405044363440678113233107768143313195\backslash\\
&  \;\;\;\;\;\;65155769860596260007646063875144448165163256825025...,
\end{align*}

\begin{align*}
A\left(  \frac{1}{3}\right)   &  =A\left(  \frac{2}{3}\right)
=2.12911826567282015384415646541541766703489387405027\backslash\\
&
\;\;\;\;\;\;\;\;\;\;\;\;\;\;\;\;\;\;\;\;\;\;79288309048724472004608367819324592064930700131124...,
\end{align*}%
\begin{align*}
A\left(  \frac{1}{5}\right)   &  =A\left(  \frac{4}{5}\right)
=3.50782830348353864334882468882849726413077940203177\backslash\\
&
\;\;\;\;\;\;\;\;\;\;\;\;\;\;\;\;\;\;\;\;\;\;20807962506019404598992541609606082898534188329417....
\end{align*}

The polynomials corresponding to%
\[%
\begin{array}
[c]{ccccc}%
x_{n}=\dfrac{x_{n-1}}{1+x_{n-1}+x_{n-1}^{2}} &  & \text{for }n\geq1, &  &
0<x_{0}<\infty
\end{array}
\]
are%
\[%
\begin{array}
[c]{ccccccc}%
P_{0}=1, &  & P_{1}=X, &  & P_{2}=-\dfrac{1}{2}+X+X^{2}, &  & P_{3}=-\dfrac
{1}{6}-\dfrac{1}{2}X+\dfrac{5}{2}X^{2}+X^{3},
\end{array}
\]%
\[%
\begin{array}
[c]{ccc}%
P_{4}=\dfrac{7}{36}-\dfrac{7}{6}X+\dfrac{3}{2}X^{2}+\dfrac{13}{3}X^{3}%
+X^{4}, &  & P_{5}=\dfrac{89}{720}-\dfrac{7}{36}X-\dfrac{17}{12}X^{2}%
+\dfrac{41}{6}X^{3}+\dfrac{77}{12}X^{4}+X^{5},
\end{array}
\]%
\[
P_{6}=-\dfrac{331}{3600}+\dfrac{197}{360}X-2X^{2}+4X^{3}+\dfrac{50}{3}%
X^{4}+\dfrac{87}{10}X^{5}+X^{6}.
\]
For example,%
\begin{align*}
B\left(  1\right)   &
=0.76799378613615405044363440678113233107768143313195\backslash\\
&  \;\;\;\;\;\;65155769860596260007646063875144448165163256825025...,
\end{align*}%
\begin{align*}
B\left(  \frac{1}{2}\right)   &  =B\left(  2\right)
=1.12911826567282015384415646541541766703489387405027\backslash\\
&
\;\;\;\;\;\;\;\;\;\;\;\;\;\;\;\;\;\;\;79288309048724472004608367819324592064930700131124...,
\end{align*}%
\begin{align*}
B\left(  \frac{1}{4}\right)   &  =B\left(  4\right)
=2.50782830348353864334882468882849726413077940203177\backslash\\
&
\;\;\;\;\;\;\;\;\;\;\;\;\;\;\;\;\;\;\;20807962506019404598992541609606082898534188329417....
\end{align*}
Our conjectured identity for $A$ and $B$ starts to emerge from the mist,
verified to $100$ decimal digits. \ 

The argument of $P_{m}$ in the power-logarithmic series is now changed from
$-\ln(n)-C$ to $\ln(n)-C$. \ (What was called $b_{1}$ in \cite{F2-abel} is
$-1$ for $I$ and $J$.) \ The polynomials corresponding to%
\[%
\begin{array}
[c]{ccccc}%
x_{n}=\dfrac{x_{n-1}}{1+x_{n-1}-x_{n-1}^{2}} &  & \text{for }n\geq1, &  &
0<x_{0}<1
\end{array}
\]
are%
\[%
\begin{array}
[c]{ccccccc}%
P_{0}=1, &  & P_{1}=X, &  & P_{2}=-\dfrac{3}{2}-X+X^{2}, &  & P_{3}=\dfrac
{2}{3}-\dfrac{7}{2}X-\dfrac{5}{2}X^{2}+X^{3},
\end{array}
\]%
\[%
\begin{array}
[c]{ccc}%
P_{4}=\dfrac{121}{36}+\dfrac{37}{6}X-\dfrac{9}{2}X^{2}-\dfrac{13}{3}%
X^{3}+X^{4}, &  & P_{5}=-\dfrac{2189}{720}+\dfrac{383}{36}X+\dfrac{239}%
{12}X^{2}-\dfrac{19}{6}X^{3}-\dfrac{77}{12}X^{4}+X^{5},
\end{array}
\]%
\[
P_{6}=-\dfrac{30811}{3600}-\dfrac{10397}{360}X+12X^{2}+43X^{3}+\dfrac{5}%
{3}X^{4}-\dfrac{87}{10}X^{5}+X^{6}.
\]
For example,%
\begin{align*}
I\left(  \frac{1}{3}\right)   &
=3.48300671252821771137045207597917266865627702999515\backslash\\
&  \;\;\;\;\;\;93109810522018534103367109602141926567367016557660...,
\end{align*}%
\begin{align*}
I\left(  \frac{1}{2}\right)   &
=1.64018851423987983185892906227384028602174585754665\backslash\\
&  \;\;\;\;\;\;98725869701043989615440421793595093368942472428035...,
\end{align*}%
\begin{align*}
I\left(  \frac{2}{3}\right)   &
=0.24236645469421107664189070448882075571078568092922\backslash\\
&  \;\;\;\;\;\;33356010869980943641381391351320915203462638811790...,
\end{align*}%
\begin{align*}
I\left(  \frac{3}{4}\right)   &
=-0.47037171508299810681801917853830768357650565329723\backslash\\
&  \;\;\;\;\;\;\;\;\;21363659129577482933553934482841917192802141645832....
\end{align*}

The polynomials corresponding to%
\[%
\begin{array}
[c]{ccccc}%
x_{n}=\dfrac{x_{n-1}\left(  1+x_{n-1}\right)  }{1+2x_{n-1}} &  & \text{for
}n\geq1, &  & 0<x_{0}<\infty
\end{array}
\]
are%
\[%
\begin{array}
[c]{ccccccc}%
P_{0}=1, &  & P_{1}=X, &  & P_{2}=-\dfrac{1}{2}-X+X^{2}, &  & P_{3}=\dfrac
{2}{3}-\dfrac{1}{2}X-\dfrac{5}{2}X^{2}+X^{3},
\end{array}
\]%
\[%
\begin{array}
[c]{ccc}%
P_{4}=-\dfrac{5}{36}+\dfrac{19}{6}X+\dfrac{3}{2}X^{2}-\dfrac{13}{3}X^{3}%
+X^{4}, &  & P_{5}=-\dfrac{749}{720}-\dfrac{139}{36}X+\dfrac{77}{12}%
X^{2}+\dfrac{41}{6}X^{3}-\dfrac{77}{12}X^{4}+X^{5},
\end{array}
\]%
\[
P_{6}=\dfrac{6389}{3600}-\dfrac{857}{360}X-18X^{2}+6X^{3}+\dfrac{50}{3}%
X^{4}-\dfrac{87}{10}X^{5}+X^{6}.
\]
For example,%
\begin{align*}
J\left(  \frac{1}{2}\right)   &
=2.48300671252821771137045207597917266865627702999515\backslash\\
&  \;\;\;\;\;\;93109810522018534103367109602141926567367016557660...,
\end{align*}%
\begin{align*}
J\left(  1\right)   &
=0.64018851423987983185892906227384028602174585754665\backslash\\
&  \;\;\;\;\;\;98725869701043989615440421793595093368942472428035...,
\end{align*}

\begin{align*}
J\left(  2\right)   &
=-0.75763354530578892335810929551117924428921431907077\backslash\\
&  \;\;\;\;\;\;\;\;\;66643989130019056358618608648679084796537361188209...,
\end{align*}%
\begin{align*}
J\left(  3\right)   &
=-1.47037171508299810681801917853830768357650565329723\backslash\\
&  \;\;\;\;\;\;\;\;\;21363659129577482933553934482841917192802141645832....
\end{align*}
Our conjectured identity for $I$ and $J$ likewise breaks the surface into the
light. \ A\ formal proof of both identities would be good to see someday.

Replacing $x$ by $1/w$, the expressions%
\[%
\begin{array}
[c]{ccccc}%
\left(  \dfrac{x}{1+x\pm x^{2}}\right)  ^{-1} &  & \text{become} &  &
w+1\pm\dfrac{1}{w}.
\end{array}
\]
It is known that a solution $\neq\mp1$ (identically) of the difference
equation%
\[%
\begin{array}
[c]{ccc}%
w(z+1)=w(z)+1\pm\dfrac{1}{w(z)}, &  & z\in\mathbb{C}%
\end{array}
\]
cannot satisfy any algebraic differential equation \cite{Ki-abel, Ta-abel,
Y1-abel, Y2-abel}. \ It is still possible that $w(z)$ could be specified via
non-algebraic means. \ Moreover, any nontrivial Abel function (e.g., $A$, $B$,
$I$ or $J$) must necessarily be hypertranscendental \cite{BB-abel, Fe-abel}.
\ This is a serious blow to our quest to recognize the associated constants.
\ We shall return to this topic in Section 4.

\section{Sums and Products}

The recurrences associated with $A$ and $B$ give sequences%
\[%
\begin{array}
[c]{ccc}%
\dfrac{1}{2},\dfrac{1}{4},\dfrac{3}{16},\dfrac{39}{256},\dfrac{8463}%
{65536},\dfrac{483008799}{4294967296},\ldots &  & \text{for }x_{0}=\dfrac
{1}{2}%
\end{array}
\]
and%
\[%
\begin{array}
[c]{ccc}%
\dfrac{1}{1},\dfrac{1}{3},\dfrac{3}{13},\dfrac{39}{217},\dfrac{8463}%
{57073},\dfrac{483008799}{3811958497},\ldots &  & \text{for }x_{0}=1
\end{array}
\]
respectively. \ Numerators always match \cite{SP-abel, A1-abel}; denominators
for $A$ are simply $2^{2^{n-1}}$. \ Denominators for $B$ are connected to
another sequence \cite{A2-abel}%
\[
\dfrac{2}{1},\dfrac{2}{1},\dfrac{4}{3},\dfrac{16}{13},\dfrac{256}{217}%
,\dfrac{65536}{57073},\dfrac{4294967296}{3811958497},\ldots
\]
defined as follows. \ Let $t_{1}=2$; for $n\geq2$, $t_{n}$ has the property
that%
\[
t_{1}+t_{2}+\cdots+t_{n}=t_{1}\,t_{2}\,\cdots\,t_{n}.
\]
Using $\{t_{n}\}$, we can circle back and characterize numerators of $A$
\&\ $B$ as denominators of $%
{\textstyle\sum_{m=1}^{n}}
t_{m}$.

Likewise, for $A$ and $B$, we have%
\[%
\begin{array}
[c]{ccc}%
\dfrac{1}{3},\dfrac{2}{9},\dfrac{14}{81},\dfrac{938}{6561},\dfrac
{5274374}{43046721},\dfrac{199225484935778}{1853020188851841},\ldots &  &
\text{for }x_{0}=\dfrac{1}{3}%
\end{array}
\]
and%
\[%
\begin{array}
[c]{ccc}%
\dfrac{1}{2},\dfrac{2}{7},\dfrac{14}{67},\dfrac{938}{5623},\dfrac
{5274374}{37772347},\dfrac{199225484935778}{1653794703916063},\ldots &  &
\text{for }x_{0}=\dfrac{1}{2}%
\end{array}
\]
respectively. \ Numerators always match \cite{A3-abel}; denominators for $A$
are simply $3^{2^{n-1}}$. \ Denominators for $B$ are connected to another
sequence \cite{A4-abel}%
\[
\dfrac{3}{1},\dfrac{3}{2},\dfrac{9}{7},\dfrac{81}{67},\dfrac{6561}%
{5623},\dfrac{43046721}{37772347},\dfrac{1853020188851841}{1653794703916063}%
,\ldots
\]
defined as before except $t_{1}=3$. \ Using $\{t_{n}\}$, we can circle back
and characterize numerators of $A$ \&\ $B$ as denominators of cumulative sums.

\section{Ratios}

The recurrences associated with $I$ and $J$ give sequences%
\[%
\begin{array}
[c]{ccc}%
\dfrac{1}{2},\dfrac{2}{5},\dfrac{10}{31},\dfrac{310}{1171},\dfrac
{363010}{1638151},\dfrac{594665194510}{3146427633211},\ldots &  & \text{for
}x_{0}=\dfrac{1}{2}%
\end{array}
\]
and%
\[%
\begin{array}
[c]{ccc}%
\dfrac{1}{1},\dfrac{2}{3},\dfrac{10}{21},\dfrac{310}{861},\dfrac
{363010}{1275141},\dfrac{594665194510}{2551762438701},\ldots &  & \text{for
}x_{0}=1
\end{array}
\]
respectively. \ Numerators always match and were conjectured by Quet
\cite{A5-abel} to follow%
\[%
\begin{array}
[c]{ccc}%
u_{n+1}=u_{n}^{2}+\dfrac{u_{n}^{3}}{u_{n-1}^{2}}-u_{n}u_{n-1}^{2}, &  &
n\geq2.
\end{array}
\]
Denominators for $I$ appear to be ratios $u_{n+1}/u_{n}$ \cite{A6-abel} and
denominators for $J$ were conjectured to follow \cite{A7-abel}%
\[%
\begin{array}
[c]{ccc}%
v_{n+1}=v_{n}^{2}+\dfrac{v_{n}^{3}}{2v_{n-1}^{2}}-\dfrac{v_{n}v_{n-1}^{2}}%
{2}, &  & n\geq2.
\end{array}
\]

Likewise, for $I$ and $J$, we have%
\[%
\begin{array}
[c]{ccc}%
\dfrac{1}{3},\dfrac{3}{11},\dfrac{33}{145},\dfrac{4785}{24721},\dfrac
{118289985}{706521601},\dfrac{83574429584465985}{568754681712768961},\ldots &
& \text{for }x_{0}=\dfrac{1}{3}%
\end{array}
\]
and%
\[%
\begin{array}
[c]{ccc}%
\dfrac{1}{2},\dfrac{3}{8},\dfrac{33}{112},\dfrac{4785}{19936},\dfrac
{118289985}{588231616},\dfrac{83574429584465985}{485180252128302976},\ldots &
& \text{for }x_{0}=\dfrac{1}{2}%
\end{array}
\]
respectively. \ Numerators always match with the same second-order
$u$-recurrence as before \cite{A8-abel}. \ Denominators for $I$ appear to be
ratios of successive numerators \cite{A9-abel} and denominators for $J$ again
follow the $v$-recurrence \cite{A0-abel}.

\section{Higher Powers}

Iterates of $y/\left(  1+y+y^{2}\right)  $ satisfy%
\begin{align*}
y_{n}  &  \sim\frac{1}{n}-\frac{\ln(n)}{n^{2}}-\frac{C}{n^{2}}+\frac
{\ln(n)^{2}}{n^{3}}-(1-2C)\frac{\ln(n)}{n^{3}}-\left(  \frac{1}{2}%
+C-C^{2}\right)  \frac{1}{n^{3}}\\
&  -\frac{\ln(n)^{3}}{n^{4}}+\left(  \frac{5}{2}-3C\right)  \frac{\ln(n)^{2}%
}{n^{4}}+\left(  \frac{1}{2}+5C-3C^{2}\right)  \frac{\ln(n)}{n^{4}}+\left(
-\frac{1}{6}+\frac{1}{2}C+\frac{5}{2}C^{2}-C^{3}\right)  \frac{1}{n^{4}}\\
&  +\frac{\ln(n)^{4}}{n^{5}}+\left(  -\frac{13}{3}+4C\right)  \frac{\ln
(n)^{3}}{n^{5}}+\left(  \frac{3}{2}-13C+6C^{2}\right)  \frac{\ln(n)^{2}}%
{n^{5}}\\
&  +\left(  \frac{7}{6}+3C-13C^{2}+4C^{3}\right)  \frac{\ln(n)}{n^{5}}+\left(
\frac{7}{36}+\frac{7}{6}C+\frac{3}{2}C^{2}-\frac{13}{3}C^{3}+C^{4}\right)
\frac{1}{n^{5}}%
\end{align*}
as $n\rightarrow\infty$ and, if $y_{0}=1/2$, then $C=1.129118265672820....$
\ (Also, if $y_{0}=1$, then $C=0.767993786136154...$, the minimum value.)
\ Iterates of $y/\left(  1+y-y^{2}\right)  $ likewise satisfy%
\begin{align*}
y_{n}  &  \sim\frac{1}{n}+\frac{\ln(n)}{n^{2}}-\frac{C}{n^{2}}+\frac
{\ln(n)^{2}}{n^{3}}-(1+2C)\frac{\ln(n)}{n^{3}}+\left(  -\frac{3}{2}%
+C+C^{2}\right)  \frac{1}{n^{3}}\\
&  +\frac{\ln(n)^{3}}{n^{4}}-\left(  \frac{5}{2}+3C\right)  \frac{\ln(n)^{2}%
}{n^{4}}+\left(  -\frac{7}{2}+5C+3C^{2}\right)  \frac{\ln(n)}{n^{4}}+\left(
\frac{2}{3}+\frac{7}{2}C-\frac{5}{2}C^{2}-C^{3}\right)  \frac{1}{n^{4}}\\
&  +\frac{\ln(n)^{4}}{n^{5}}-\left(  \frac{13}{3}+4C\right)  \frac{\ln(n)^{3}%
}{n^{5}}+\left(  -\frac{9}{2}+13C+6C^{2}\right)  \frac{\ln(n)^{2}}{n^{5}}\\
&  +\left(  \frac{37}{6}+9C-13C^{2}-4C^{3}\right)  \frac{\ln(n)}{n^{5}%
}+\left(  \frac{121}{36}-\frac{37}{6}C-\frac{9}{2}C^{2}+\frac{13}{3}%
C^{3}+C^{4}\right)  \frac{1}{n^{5}}%
\end{align*}
and, if $y_{0}=1/2$, then $C=1.640188514239879....$ \ Another way of defining
$C$ involves $x_{n}=1/y_{n}$. \ That is, the recurrence%
\[
x_{n+1}=x_{n}+1\pm\frac{1}{x_{n}}%
\]
has asymptotics%
\[
\lim_{n\rightarrow\infty}\left(  x_{n}-n\mp\ln(n)\right)  =C
\]
or, encompassing greater detail \cite{P1-abel, P2-abel},%
\[
x_{n}\sim n+\ln(n)+C+\frac{\ln(n)}{n}+\left(  \frac{1}{2}+C\right)  \frac
{1}{n}%
\]
and%
\[
x_{n}\sim n-\ln(n)+C+\frac{\ln(n)}{n}+\left(  \frac{3}{2}-C\right)  \frac
{1}{n}%
\]
respectively. \ Given $x_{0}=2$, the values $1.129118265672820...$ and
$1.640188514239879...$ \ for $C$ again apply. \ A\ division-based technique
for calculating more terms in the series expansion of $x_{n}$ is not known;
finding the reciprocal of a power-logarithmic series seems to be generally difficult.

In contrast, iterates of $y/\left(  1+y+y^{3}\right)  $ satisfy%
\begin{align*}
y_{n}  &  \sim\frac{1}{n}-\frac{C}{n^{2}}+\left(  1+C^{2}\right)  \frac
{1}{n^{3}}+\left(  \frac{1}{2}-3C-C^{3}\right)  \frac{1}{n^{4}}+\left(
\frac{11}{6}-2C+6C^{2}+C^{4}\right)  \frac{1}{n^{5}}\\
&  +\left(  \frac{9}{4}-\frac{55}{6}C+5C^{2}-10C^{3}-C^{5}\right)  \frac
{1}{n^{6}}+\left(  \frac{299}{60}-\frac{27}{2}C+\frac{55}{2}C^{2}%
-10C^{3}+15C^{4}+C^{6}\right)  \frac{1}{n^{7}}%
\end{align*}
and, if $y_{0}=1/2$, then $C=2.598786855824871....$ \ (Also, if $y_{0}%
=2^{-1/3}$, then $C=2.286858220891602...$, the minimum value.) \ Iterates of
$y/\left(  1+y-y^{3}\right)  $ satisfy%
\begin{align*}
y_{n}  &  \sim\frac{1}{n}-\frac{C}{n^{2}}-\left(  1-C^{2}\right)  \frac
{1}{n^{3}}-\left(  \frac{1}{2}-3C+C^{3}\right)  \frac{1}{n^{4}}+\left(
\frac{3}{2}+2C-6C^{2}+C^{4}\right)  \frac{1}{n^{5}}\\
&  +\left(  \frac{9}{4}-\frac{15}{2}C-5C^{2}+10C^{3}-C^{5}\right)  \frac
{1}{n^{6}}-\left(  \frac{27}{20}+\frac{27}{2}C-\frac{45}{2}C^{2}%
-10C^{3}+15C^{4}-C^{6}\right)  \frac{1}{n^{7}}%
\end{align*}
and, if $y_{0}=1/2$, then $C=1.290937947423058....$ \ No logarithmic terms
appear in either expansion. \ (Reason:\ what was called $b_{1}$ in
\cite{F2-abel} is $0$ here.) \ Finding the reciprocal of an ordinary power
series is more easily done than for a power-logarithmic series. \ The
recurrence%
\[
x_{n+1}=x_{n}+1\pm\frac{1}{x_{n}^{2}}%
\]
has asymptotics
\begin{align*}
x_{n}  &  \sim n+C-\frac{1}{n}-\left(  \frac{1}{2}-C\right)  \frac{1}{n^{2}%
}+\left(  -\frac{5}{6}+C-C^{2}\right)  \frac{1}{n^{3}}+\left(  -\frac{5}%
{4}+\frac{5}{2}C-\frac{3}{2}C^{2}+C^{3}\right)  \frac{1}{n^{4}}\\
&  +\left(  -\frac{31}{15}+5C-5C^{2}+2C^{3}-C^{4}\right)  \frac{1}{n^{5}%
}+\left(  -\frac{11}{3}+\frac{31}{3}C-\frac{25}{2}C^{2}+\frac{25}{3}%
C^{3}-\frac{5}{2}C^{4}+C^{5}\right)  \frac{1}{n^{6}}%
\end{align*}
and
\begin{align*}
x_{n}  &  \sim n+C+\frac{1}{n}+\left(  \frac{1}{2}-C\right)  \frac{1}{n^{2}%
}+\left(  -\frac{1}{2}-C+C^{2}\right)  \frac{1}{n^{3}}+\left(  -\frac{5}%
{4}+\frac{3}{2}C+\frac{3}{2}C^{2}-C^{3}\right)  \frac{1}{n^{4}}\\
&  +\left(  -\frac{2}{5}+5C-3C^{2}-2C^{3}+C^{4}\right)  \frac{1}{n^{5}%
}+\left(  \frac{5}{2}+2C-\frac{25}{2}C^{2}+5C^{3}+\frac{5}{2}C^{4}%
-C^{5}\right)  \frac{1}{n^{6}}%
\end{align*}
respectively. \ Given $x_{0}=2$, the values $2.598786855824871...$ and
$1.290937947423058...$ \ for $C$ again apply.

Ordinary power series arise for%
\[
x_{n+1}=x_{n}+1\pm\frac{1}{x_{n}^{\ell}},
\]
given any integer $\ell\geq2$. \ One might expect that, for such iterations,
determing a closed-form expression for $C$ might be achievable. \ This belief
is challenged by a result \cite{Y3-abel} that the difference equation%
\[%
\begin{array}
[c]{ccc}%
w(z+1)=w(z)+1\pm\dfrac{1}{w(z)^{\ell}}, &  & z\in\mathbb{C}%
\end{array}
\]
possesses no differentially algebraic solutions. \ Thus the striking
distinction between the cases $\ell=1$ and $\ell\geq2$ is immaterial when
attempting to understand $C$. \ Again, no nontrivial Abel function is
differentially algebraic \cite{BB-abel, Fe-abel}. \ While capturing $C$ still
remains a potent aspiration, the task\ seems increasingly out of reach.%

\begin{figure}
[ptb]
\begin{center}
\includegraphics[
height=2.8504in,
width=2.3471in
]%
{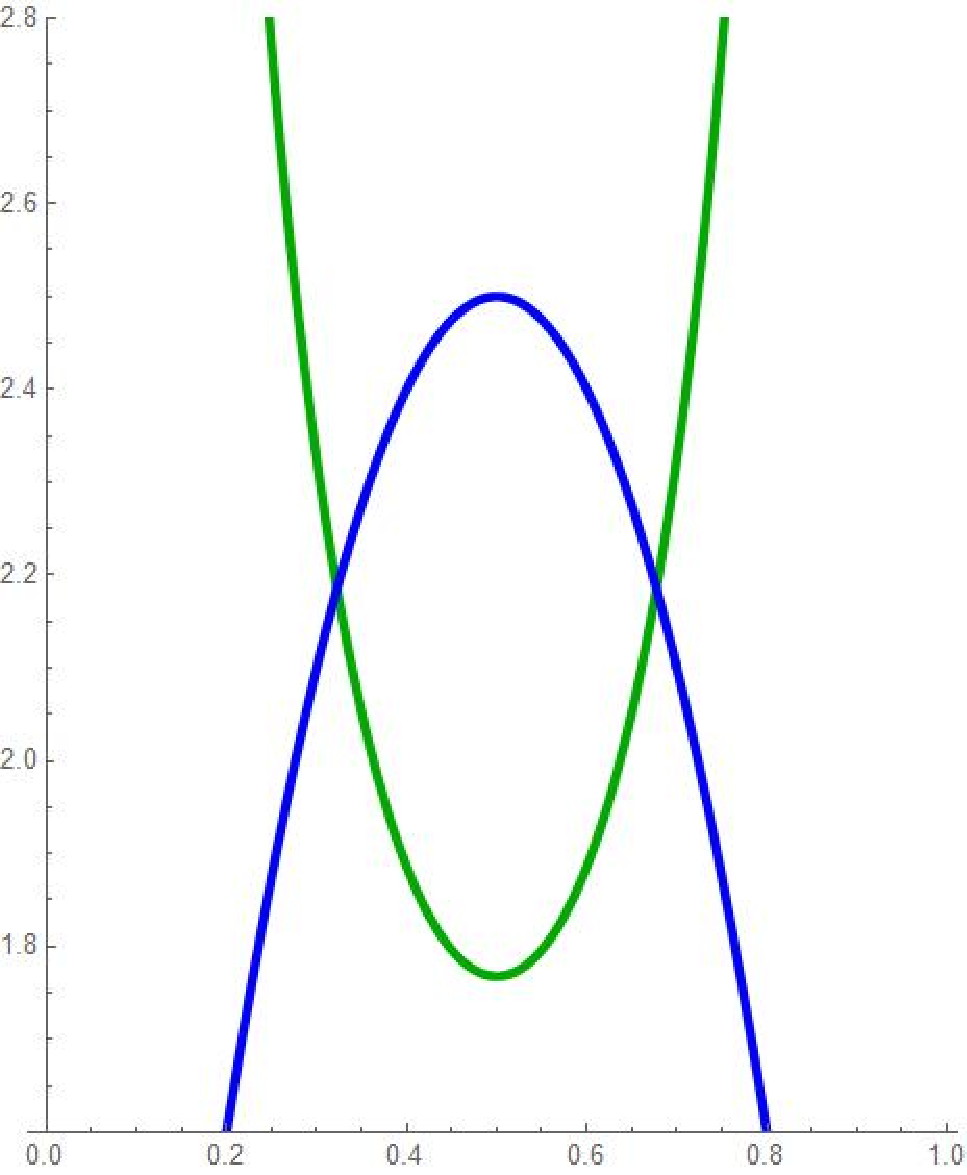}%
\caption{Blue curve is $10x(1-x)$, scaled for visibility, with maximum at
$x=1/2$. \ Green curve is $A(x)$, with minimum at $x=1/2$. \ Axes are
compatible.}%
\end{center}
\end{figure}
\begin{figure}
[ptb]
\begin{center}
\includegraphics[
height=2.8504in,
width=2.725in
]%
{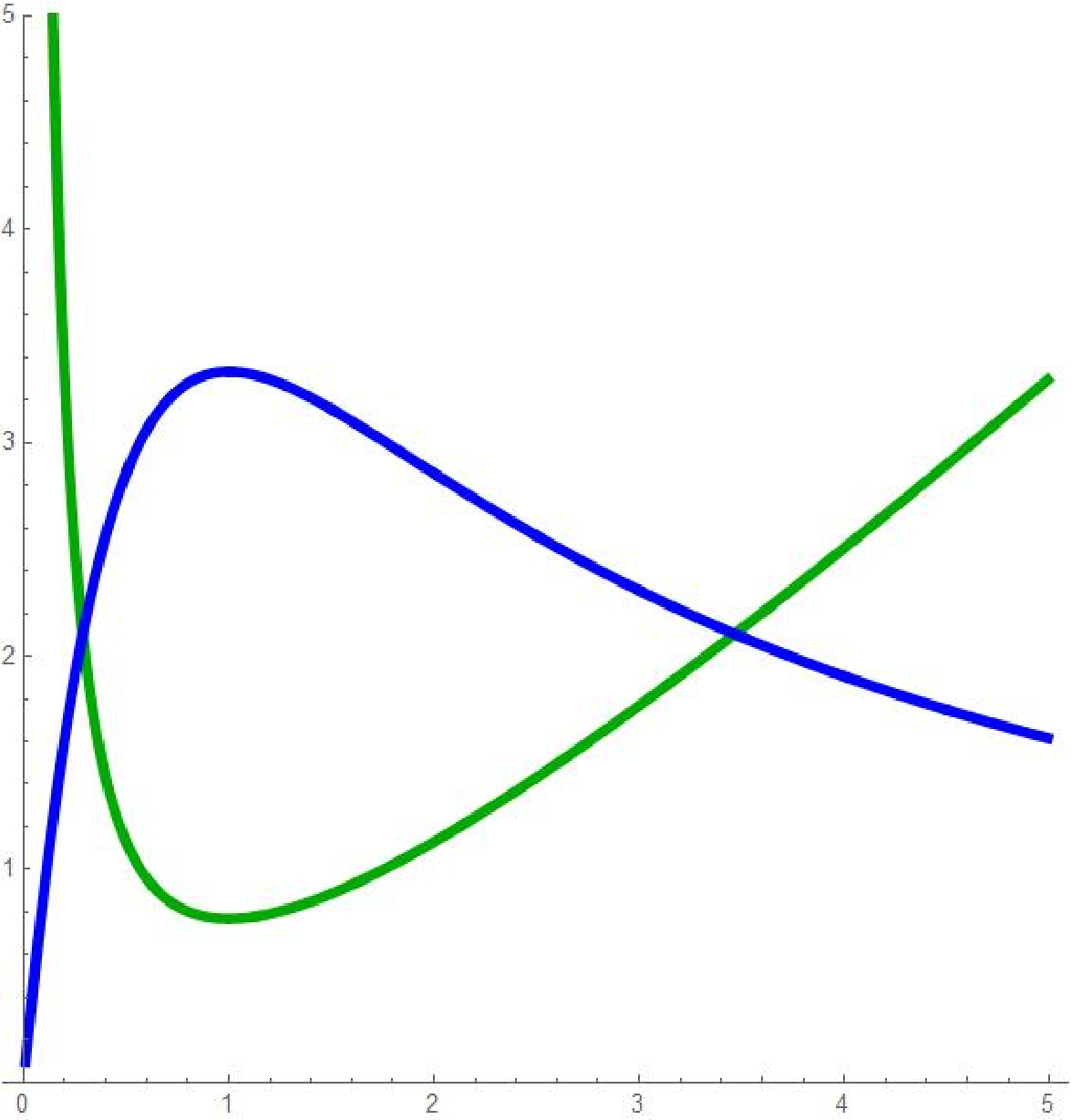}%
\caption{Blue curve is $10x/(1+x+x^{2})$, scaled for visibility, with maximum
at $x=1$. \ Green curve is $B(x)$, with minimum at $x=1$. \ Axes are
compatible.}%
\end{center}
\end{figure}
\begin{figure}
[ptb]
\begin{center}
\includegraphics[
height=2.8504in,
width=2.8253in
]%
{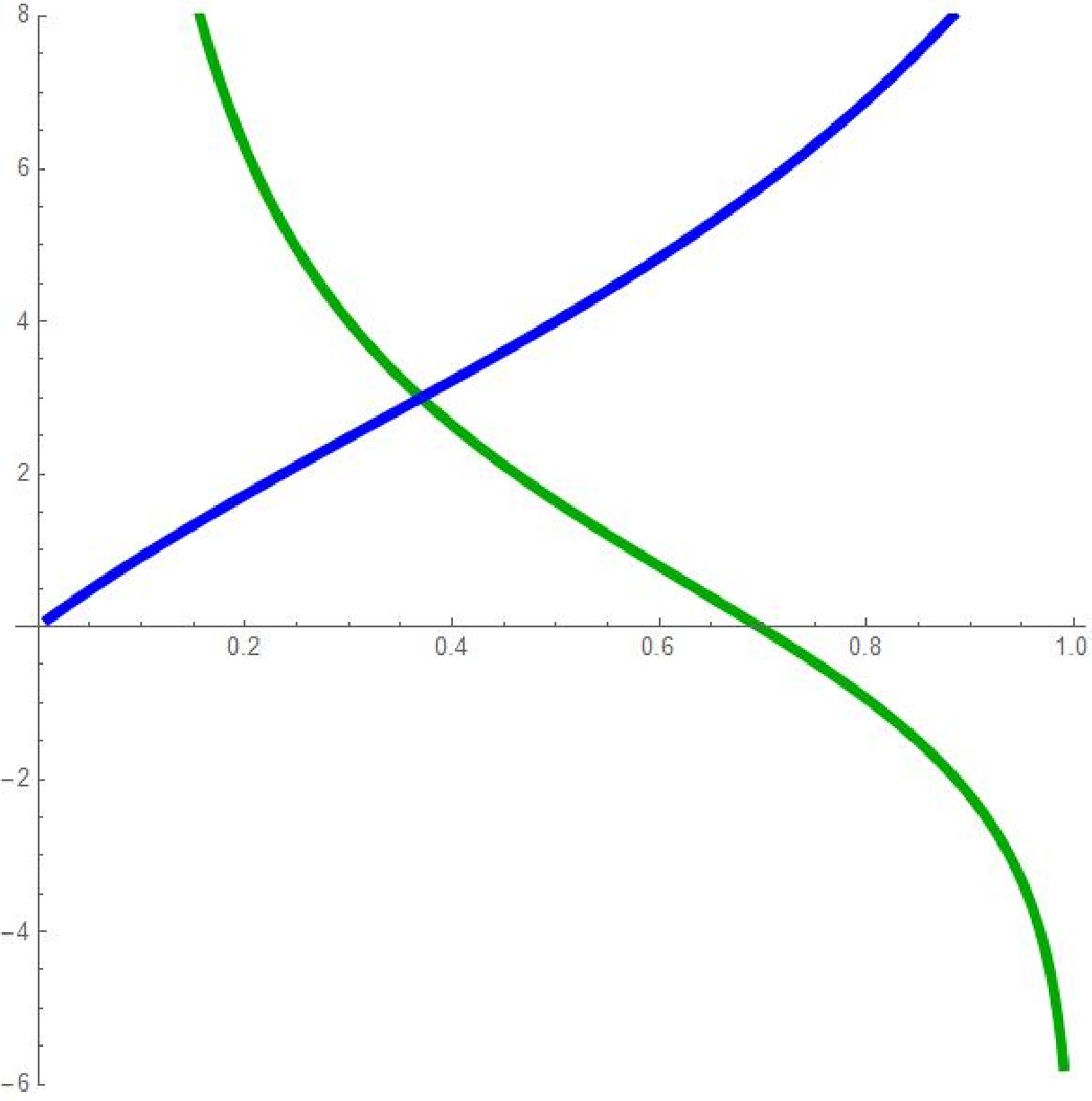}%
\caption{Blue curve is $10x/(1+x-x^{2})$, scaled for visibility, with
inflection point at $\varphi^{1/3}-\varphi^{-1/3}=0.322$. \ Green curve is
$I(x)$, with inflection point at $\psi=0.629$. \ Distance between vertical
axis notches is ten times that for horizontal. }%
\end{center}
\end{figure}
\begin{figure}
[ptb]
\begin{center}
\includegraphics[
height=2.8504in,
width=2.8167in
]%
{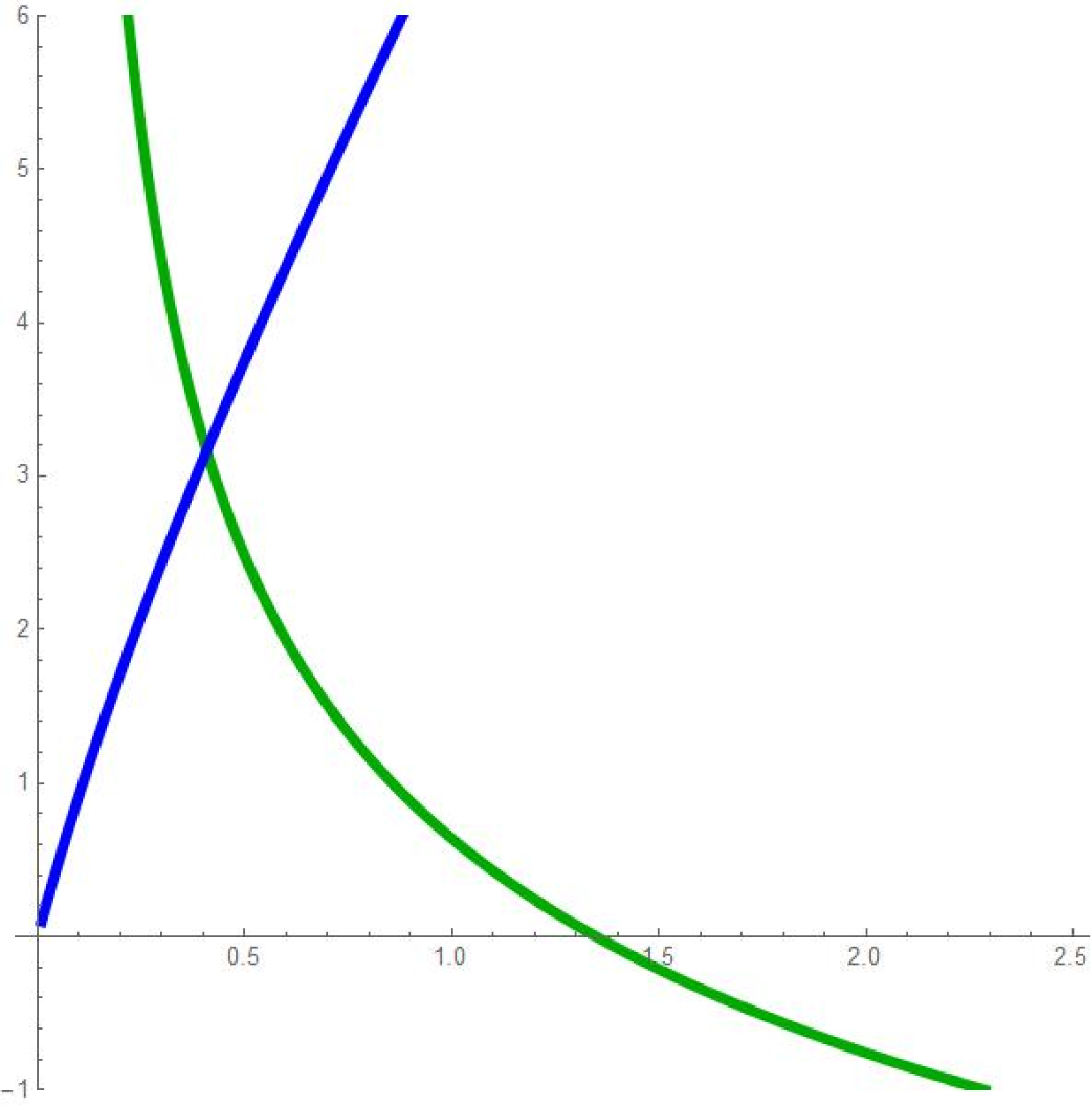}%
\caption{Blue curve is $10x(1+x)/(1+2x)$, scaled for visibility. Green curve
is $J(x)$. Distance between vertical axis notches is two times that for
horizontal.}%
\end{center}
\end{figure}

\section{Corrigendum}

The interrelation between $A$ and $B$ turns out to be frivolous. \ Starting
with the $B$-recurrence%
\[%
\begin{array}
[c]{ccc}%
x_{n+1}=\dfrac{x_{n}}{1+x_{n}+x_{n}^{2}}, &  & \xi_{n}=\dfrac{1}{x_{n}}%
\end{array}
\]
and reciprocating, we obtain%
\[%
\begin{array}
[c]{ccc}%
\xi_{n+1}=\xi_{n}+1+\dfrac{1}{\xi_{n}}, &  & \tilde{x}_{n}=\dfrac{1}{1+\xi
_{n}}.
\end{array}
\]
Because $\xi_{n}=1/\tilde{x}_{n}-1$, it follows that
\[
\dfrac{1-\tilde{x}_{n}}{\tilde{x}_{n}}+1+\dfrac{\tilde{x}_{n}}{1-\tilde{x}%
_{n}}=\dfrac{\left(  1-\tilde{x}_{n}\right)  ^{2}+\tilde{x}_{n}\left(
1-\tilde{x}_{n}\right)  +\tilde{x}_{n}^{2}}{\tilde{x}_{n}\left(  1-\tilde
{x}_{n}\right)  }=\dfrac{1}{\tilde{x}_{n}\left(  1-\tilde{x}_{n}\right)
}-1=\xi_{n+1}%
\]
if and only if the $A$-recurrence holds:%
\[
\tilde{x}_{n+1}=\tilde{x}_{n}\left(  1-\tilde{x}_{n}\right)
\]
which demonstrates that one is simply a reparametrization of the other. \ The
same is true for $I$ and $J$. \ Starting with the $J$-recurrence%
\[%
\begin{array}
[c]{ccc}%
x_{n+1}=\dfrac{x_{n}(1+x_{n})}{1+2x_{n}}, &  & \xi_{n}=\dfrac{1}{x_{n}}%
\end{array}
\]
and reciprocating, we obtain%
\[%
\begin{array}
[c]{ccc}%
\xi_{n+1}=\dfrac{\xi_{n}\left(  2+\xi_{n}\right)  }{1+\xi_{n}}, &  & \tilde
{x}_{n}=\dfrac{1}{1+\xi_{n}}.
\end{array}
\]
It follows that%
\[
\dfrac{\left(  \dfrac{1}{\tilde{x}_{n}}-1\right)  \left(  2+\dfrac{1}%
{\tilde{x}_{n}}-1\right)  }{1+\dfrac{1}{\tilde{x}_{n}}-1}=\dfrac{1-\tilde
{x}_{n}^{2}}{\tilde{x}_{n}}=\dfrac{1+\tilde{x}_{n}-\tilde{x}_{n}^{2}}%
{\tilde{x}_{n}}-1=\xi_{n+1}%
\]
if and only if the $I$-recurrence holds:%
\[
\tilde{x}_{n+1}=\dfrac{\tilde{x}_{n}}{1+\tilde{x}_{n}-\tilde{x}_{n}^{2}}.
\]
Although these late realizations sadly eliminate the motivating purpose of my
paper, I\ believe that there is still some value in keeping it visible.

\section{Acknowledgements}

I am grateful to Daniel Lichtblau at Wolfram Research for kindly answering my
questions, e.g., about generalizing my original Mathematica code to arbitrary
$k$. \ The creators of Mathematica and of the On-Line Encyclopedia of Integer
Sequences earn my gratitude every day:\ this paper could not have otherwise
been written. \ An interactive computational notebook is available
\cite{F3-abel} which might be helpful to interested readers.

\end{document}